\documentclass[12pt]{amsart}

\usepackage[all]{xy}
\usepackage{mathrsfs}
\usepackage{amsmath}
\usepackage{latexsym}
\usepackage{amssymb}
\usepackage{amsfonts}
\usepackage{longtable}
\usepackage{colortbl}
\include{PDF}
\usepackage{graphics}

\renewcommand{\proof}{\par\noindent{\it Proof.\ \ }}
\def\qed{\ifmmode\square\else\nolinebreak\hfill
$\Box$\fi\par\vskip12pt}

\def\ov{\overline} 
\def\l{\langle} \def\r{\rangle} 
\def\div{\,\big|\,} 

\def\mod{{\rm mod~}}

\def\FF{\mathbb F} \def\ZZ{{\mathbb Z}}
\def\BB{{\mathcal B}}

\def\calB{{\mathcal B}}

\def\calO{{\mathcal O}}

\def\ZZ{{\mathbb Z}}

\def\mod{{\rm mod~}}

\def\K{{\bf K}}

\def\Aut{{\rm Aut}}  \def\Out{{\rm Out}}

\def\Cay{{\rm Cay}}

\def\D{{\rm D}} 
\def\S{{\rm S}} 
 \def\M{{\rm M}}
\def\soc{{\rm soc}} 
\def\C{{\bf C}}\def\Z{{\bf Z}} 

\def\Ga{\Gamma}
\def\Ome{\Omega}
\def\Sig{\Sigma}
\def\Del{\Delta}

\def\a{\alpha} \def\b{\beta} \def\d{\delta} \def\s{\sigma}
 
\def\o{\omega}

\def\GammaL{{\rm \Gamma L}}

\def\PGammaL{{\rm P\Gamma L}} 
\def\A{{\rm A}}
\def\Sym{{\rm Sym}}
\def\Alt{{\rm Alt}}

\def\PSL{{\rm PSL}}\def\PGL{{\rm PGL}}
\def\GL{{\rm GL}} 
\def\AGL{{\rm AGL}}

  \def\D{{\rm D}}

\newtheorem{theorem}{Theorem}[section]%
\newtheorem{lemma}[theorem]{Lemma}%
\newtheorem{proposition}[theorem]{Proposition}%
\newtheorem{example}[theorem]{Example}%
\newtheorem{construction}[theorem]{Construction}%
\newtheorem{hypothesis}[theorem]{Hypothesis}%

\input amssym.def
\input amssym.tex

\begin{document}

\title[Finite d-groups]
{Finite permutation groups containing a regular dihedral subgroup}
\thanks{1991 MR Subject Classification 20B15, 20B30, 05C25.}
\thanks{This work is partially supported by  NSFC (No. 61771019)}

\author{Shujiao Song}
\address{School of Mathematics\\
Yantai University\\
China}
\email{}

\date\today

\maketitle

\date\today

\begin{abstract}
A characterization is given of finite permutation groups which contain a regular dihedral subgroup.
\end{abstract}

\section{Introduction}

A permutation group $G\le\Sym(\Ome)$ is called a {\it c-group} or a {\it d-group} if it contains a regular subgroup which is cyclic or dihedral, respectively.
The study of these two classes of permutation groups has a long history, dating back to Burnside \cite{Burnside} who proved that a primitive c-group of composite prime-power degree is necessarily 2-transitive.
Schur (1933, see \cite{Wielandt}) completed Burnside's work by proving that a primitive c-group of composite degree is 2-transitive.
Based on the Classification of Finite Simple Groups, primitive simple c-groups were listed in \cite{Feit,Gorenstein,Kantor}.
More recently, Jones \cite{Jones} and Li \cite{Abel-B-gps} completely determined primitive permutation c-groups, independently.
Then in 2010, Li and Praeger \cite{c-gps} studied general c-groups.

Wielandt proved that primitive d-groups are 2-transitive, refer to \cite{Wielandt}, extended Schur's result regarding c-groups to d-groups.
Then in \cite{et-Cay,Pan-2,Song-Li}, quasiprimitive d-groups and bi-quasiprimitive d-groups are classified.
The main theorem of this paper gives a recursive characterization of $d$-groups.

The recent research of c-groups and d-groups is often associated with the study of symmetrical graphs (directed or undirected).
A {\it circulant} is a Cayley graph of a cyclic group, and a {\it dihedrant} is a Cayley graph of a dihedral group.
An {\it arc} of a graph is an ordered pair of adjacent vertices, and a {\it $2$-arc} is a triple of vertices one of which is adjacent to the other two.
A graph is called {\it arc-transitive} or {\it $2$-arc-transitive} if all of its arcs or 2-arcs, respectively, are equivalent under automorphisms.
For arc-transitive circulants, recursive characterizations are obtained in \cite{Kovacs,Li-circulant}, and an explicit classification is given in \cite{LXZ}.
For dihedrants, special classes of arc-transitive dihedrants are characterized in \cite{Du-MM,KMM,Marusic,Pan}.
However, the problem of {\bf classifying arc-transitive dihedrants} is still an unsettled open problem.
Theorem~\ref{thm-1}, we believe, would play an important role in the study of arc-transitive dihedrants.

To state the theorem, we need introduce a few concepts regarding permutation groups.
Let $G$ be a transitive permutation group on a set $\Ome$.
For a block system $\calB$ and a block $B\in\calB$, the group $G$ induces a transitive permutation group $G^\calB$ on $\calB$, and the stabilizer $G_B$ induces a transitive group $G_B^B$ on $B$.
Then $G^\calB\cong G/G_{(\calB)}$, and $G_B^B\cong G_B/G_{(B)}$, where $G_{(\calB)}$ is the kernel of $G$ on $\calB$, and $G_{(B)}$ fixes $B$ pointwise.
A block system $\calB$ is called {\it maximal} if the induced action $G^\calB$ is primitive, and called
{\it minimal} if, for a block $B\in\calB$, the induced action $G_B^B$ is primitive.

An {\it orbital graph} of $G$ is a (di)graph with vertex set $\Ome$ and edge set $(\o_1,\o_2)^G$, where $\o_1,\o_2\in\Ome$.
An orbital graph is obviously arc-transitive.
For a block system $\calB$, an orbital graph $\Ga$ has a {\it quotient graph} $\Ga_\calB$, which is an orbital graph of the induced permutation group $G^\calB$.
As usual, $\K_n$ denotes a complete graph of order $n$, and $\ov\K_n$ denotes the complement.
For a graph $\Sigma$, denote by $\Sigma[\ov\K_n]$ the {\it lexicographic product} of $\ov\K_n$ by $\Sigma$.

Since primitive d-groups are 2-transitive and completely known by Theorem~\ref{qp-d-gps}, the main result of this paper focuses on the imprimitive case.

\begin{theorem}\label{thm-1}
Let $G$ be an imprimitive d-group of degree $n$ and $\calB$ a minimal block system, and let $K=G_{(\calB)}$ and $B\in\calB$.
Then each connected orbital graph is undirected, $G^\calB$ is a c-group or a d-group, and further, one of the following holds:
\begin{itemize}
 \item[(1)] $G\cong G^\calB$, and $G$ is also a c-group of degree $n/2$, where $|\calB|={1\over2}|\Ome|$;

 \item[(2)] there exists an orbital graph for $G$ of the form $\Ga_\calB[\ov\K_b]$, $K^B$ is primitive and $K\not\cong K^B$, where $b=|B|$ and $\Ga_\calB$ is an orbital graph of $G^\calB$;

\item[(3)] $G=(K\times H).\calO$, where $K\cong K^B$ is primitive, either of $K.\calO$ and $H.\calO$ is a c-group or a d-group, and $\calO\leq\Out(K)$;

\item[(4)] $\ZZ_2=K=\Z(G)$, and $G^\calB\cong G/\Z(G)$ is a d-group of degree $n/2$, and $K.G^\calB$ is a non-split extension;

\item[(5)] $K=\ZZ_p$ with $p$ being an odd prime, and $G=\C_G(K).\ZZ_{2\ell}$, where $2\ell\div(p-1)$ and $\ell$ is odd, and $K.(\C_G(K))^\calB$ is a non-split extension;

\item[(6)] $|B|=4$, $R\cap K=\ZZ_2$, and $G=2^d.G^\calB$, where $G^\calB$ is a c-group of degree $n/4$.
\end{itemize}
\end{theorem}

In subsequent work, it will be shown that an arc-transitive dihedrant either can be explicitly reconstructed from a smaller circulant or dihedrant, or is an orbital graph of a group in part~(6).

\vskip0.08in
{\bf Remarks on Theorem~\ref{thm-1}:}
\begin{itemize}
\item[(a)] For the case in Theorem~\ref{thm-1}\,(1), the d-group $G$ of degree $n$ has a faithful permutation representation to be a c-group of degree $n/2$, which is better understood than the general case.

\item[(b)] In the view point of arc-transitive dihedrants, the case in Theorem~\ref{thm-1}\,(2) is completely reduced to smaller dihedrants.

\item[(c)] In case (3), as $\calO$ is relatively very small, $G$ is essentially determined by a direct product of smaller c-groups or d-groups.

\item[(d)] Due to the classification of primitive c-groups and d-groups, the case in Theorem~\ref{thm-1}\,(4) is completely reduced to a smaller case.

\item[(e)] The case in Theorem~\ref{thm-1}\,(3) is reduced to c-groups.

\item[(f)] The case in Theorem~\ref{thm-1}\,(5) is not well-understood as other cases, which seems challenge to obtain a nice description for this case.
\end{itemize}

The notation used in this paper is standard, see \cite{Biggs,DM-book} for example.

\section{Preliminaries}

We begin with a technical lemma for general permutation groups.

\begin{lemma}\label{key-lem}
Let $G$ be a transitive permutation group on $\Ome$ which has a connected orbital graph $\Ga$.
Let $\calB$ be a minimal block system for $G$ on $\Ome$, $B\in\calB$, and let $K=G_{(\calB)}$.
Then one of the following holds:
\begin{itemize}
\item[(i)] $G\cong G^\calB$,
\item[(ii)] $K\cong K^B$ is a transitive permutation group,
\item[(iii)] $\Ga=\Ga_\calB[\ov\K_b]$, where $b=|B|$.
\end{itemize}
\end{lemma}
\proof
If $K=1$, then $G^\calB\cong G/K\cong G$, as in part~(i).

Assume that $K\not=1$.
Then $K^B\not=1$ since the actions of $K$ on all of the blocks in $\calB$ are equivalent.
Since $K^B\lhd G_B^B$ and $G_B^B$ is primitive, $K^B$ is transitive.
If $K\cong K^B$, then $K$ has a faithful transitive representation on $B$, as in part~(ii).

Finally, assume that $K\not\cong K^B$.
Then the kernel $K_{(B)}\not=1$, and $K_{(B)}$ is non-trivial on some block $B'\in\calB$.
Since $\Ga$ is connected, there is a path $\a_0,\a_1,\dots,\a_\ell$ of $\Ga$, where
$\a_0\in B$, $\a_i\in B_i$, and $\a_\ell\in B'$, and there exists $1\leq k\leq \ell$ such that
$K_{(B)}$ is trivial on $B_i$ for $i\leq k-1$ and non-trivial on $B_k$.
Then $K_{(B_{k-1})}=K_{(B)}$, and $1\not=K_{(B)}^{B_k}=K_{(B_{k-1})}^{B_k}\lhd K^{B_k}$.
%{\bf Since $K^{B_k}$ is transitive, hence $K_{(B_{k-1})}^{B_k}$ is half-transitive. Suppose its orbit on $B_{k}$ has size $b$. Then $b\mid |B|$, and so $[B_{k-1},B_K]\cong{ |B|\over b} \K_{b,b}$}

Since $K^{B_k}$ is primitive, $K_{(B)}^{B_k}$ is transitive, and hence the induced
subgraph on $[B_{k-1},B_k]$ is isomorphic to $\K_{b,b}$ where $b=|B_{k-1}|=|B_k|=|B|$.
Since $\Ga$ is $G$-arc-transitive, it follows that $\Ga=\Ga_\calB[\ov\K_b]$.
\qed

In general, an orbital graph is not necessarily undirected.
However, the next lemma shows that connected orbital graphs of a d-group are always undirected.

\begin{lemma}\label{orbital-d-gps}
Let $G$ be a d-group on $\Ome$, and $\Ga$ an orbital graph of $G$.
Then either
\begin{itemize}
\item[(i)] $\Ga$ is an undirected arc-transitive graph, or
\item[(ii)] $\Ga$ is disconnected, and each component is a circulant.
\end{itemize}
Furthermore, is $\Ga$ is connected, then $\Ga$ is undirected.
\end{lemma}
\proof
Let $G$ be a d-group on $\Ome$ and $\Ga$ a connected  orbital graph of $G$.
Then $G$ has a dihedral regular subgroup $D$, and $\Ga$ is a Cayley graph of $D$, namely,
$\Ga=\Cay(D,S)$ for some subset $S\subset D$.
Let $D=\l a\r{:}\l z\r$, where $|z|=2$ and $a^z=a^{-1}$.

Assume first that $\l S\r\leq\l a\r$.
Then $\Ga$ is disconnected, and each component is a Cayley graph of a subgroup of $\l a\r$, and so it is a circulant, as in part~(ii).

Now assume that $\l S\r$ is not a subgroup of $\l a\r$.
Then $S$ contains an involution $g$ in $D\setminus\l a\r$.
The permutation of $\Ome$ induced by the right multiplication of $g$ is an element of $G$, and interchanges the arcs $(1_D,g)$ and $(g,1_D)$, where $1_D$ is the identity of the group $D$.
It follows that $\Ga=\Cay(D,S)$ is undirected.
\qed

Next, we study d-groups.
%First we consider primitive d-groups, and primitive c-groups which we need.
By the well-known Schur's theorem (\cite[Theorem 25.3]{Wielandt}) and Wielandt's theorem (\cite[Theorem 25.6]{Wielandt}), we have the following result.

\begin{theorem}\label{mini-block}{\rm(Schur, Wielandt)}
If $G$ be a finite primitive c-group or d-group on a set $\Ome$, then either
\begin{itemize}
\item[(1)] $|\Ome|=p$ with $p$ prime and $G\le\AGL(1,p)$, or
\item[(2)] $G$ is $2$-transitive.
\end{itemize}
\end{theorem}

A transitive permutation group is called {\it quasiprimitive} if each of its non-trivial normal subgroup is transitive.
Clearly, a primitive permutation group is quasiprimitive, but a quasiprimitive group is not necessarily primitive.
However, for c-groups and d-groups, the quasiprimitivity and the primitivity are equivalent, and moreover,
quasiprimitive c-groups and d-groups are classified, stated in the next two theorems.
We remark that, although the proof of Theorem~\ref{mini-block}, consisting of both Schur's work and Wielandt's work, is `elementary' and elegant, the following two theorems indeed depend on the Classification of Finite Simple Groups.

\begin{theorem}\label{qp-c-gps}{\rm (see \cite{Jones,Abel-B-gps})}
Let $X\le\Sym(\Ome)$ be a quasiprimitive c-group of degree $n$. Then
$X$ is primitive on $\Ome$, and either $n=p$ is prime and
$X\le\AGL(1,p)$, or $X$ is $2$-transitive, listed in the following
table.
\end{theorem}
\[\begin{array}{|llll|}\hline
X & X_\o & n & \mbox{condition}\\ \hline

\A_n & \A_{n-1} & & n\mbox{ is odd, }n\geq 4 \\

\S_n & \S_{n-1} & & n\geq 4\\

\PGL(d,q).o & P_1 \mbox{(parabolic)} & {q^d-1\over q-1} & o\leq \PGammaL(d,q)/\PGL(d,q)\\

%\PGL(d,r^d),\ \PGL(d,r^d).d {\bf (delete this column)}& P_1 \mbox{(parabolic)} & {r^{d^2}-1\over r^d-1} & \\

\PSL(2,11) & \A_5 & 11 &\\

\M_{11} & \M_{10} & 11 &\\

\M_{23} & \M_{22} & 23 &\\ \hline
\end{array}\]
%\centerline{{\bf Table~\ref{qp-c-gps}}}

\vskip0.1in

\begin{theorem}\label{qp-d-gps}
If $X\le\Sym(\Ome)$ is a quasiprimitive d-group with a regular
dihedral subgroup $G$, then $X$ is $2$-transitive on $\Ome$ and
$(X,G,X_\o)$ is one of the triples in the following table.
\end{theorem}

\[\begin{array}{|llll|} \hline
X & G & X_\o & \mbox{conditions} \\ \hline

\A_4 & \D_4 & \ZZ_3 & \\

\S_4 & \D_4 & \S_3 & \\

\AGL(3,2) & \D_8 & \GL(3,2) & \\ \hline

\AGL(4,2) & \D_{16} & \GL(4,2) & \\

2^4\rtimes \A_7 & \D_{16} & \A_7 & \\

2^4\rtimes \S_6 & \D_{16} & \S_6 & \\

2^4\rtimes \A_6 & \D_{16} & \A_6 & \\

2^4\rtimes \S_5 & \D_{16} & \S_5 & \\

2^4\rtimes \GammaL(2,4) & \D_{16} & \GammaL(2,4) & \\ \hline

\M_{12} & \D_{12} & \M_{11} &\\

\M_{22}.2 &\D_{22} & \PSL(3,4).2 & \\

\M_{24} & \D_{24} &\M_{23}& \\

\S_{2m} & \D_{2m} & \S_{2m-1} & \\

\A_{4m} & \D_{4m} & \A_{4m-1} & \\

\PSL(2,r^f).o & \D_{r^f+1} & \ZZ_r^f\rtimes \ZZ_{\frac{r^f-1}2}.o &
r^f\equiv 3\pmod 4, o\le\ZZ_2\times\ZZ_f \\

\PGL(2,r^f).\ZZ_e & \D_{r^f+1} & \ZZ_r^f\rtimes \ZZ_{r^f-1}\rtimes \ZZ_e & r^f\equiv
1\pmod 4, e\div f \\ \hline
\end{array}\]
\vskip0.25in

We now present a method for constructing d-groups by wreath product.

\begin{construction}\label{cons-1}
{\rm
Let $H$ be a c-group on $\Del$ of degree $k$, and label $\Ome=\{1,2,\dots,n\}$ ($n$ is even) , and $G$ contains a regular dihedral group on $\Ome$, say $D=D_{2n}.$
Let
\[X=H\wr G=H^n\rtimes G=(H_1\times H_2\times\dots\times H_n)\rtimes G.\]
Let $\Sig=\{(\d_1,\d_2,\dots,\d_n)\mid \d_i\in\Del\}$, and let
\[L=((H_1)_{\d_1}\times H_2\times\dots\times H_n)\rtimes D.\]
Then $L$ is a stabilizer of the group $X$ acting on the set $\Sig$.
}
\end{construction}

\begin{lemma}\label{lexi-d-gps}
A group $X$ constructed in Construction~$\ref{cons-1}$ is a d-group on $\Sig$ of degree $kn$.
\end{lemma}
\proof
 Let $n=2m$.
Without loss of generality, write $D=\l a\r\rtimes \l b\r$, where
$a=(12\dots m)(m+1,m+2,\dots,2m)$ and $b=(1,2m)(2,2m-1)\dots(m,m+1)$.
Let $x=(h_1,\dots,h_n)$ such that $h_1=1$, $h_{m+1}=-1$, and other $h_i$ equal 0, and let
\[y=xa=(1,0,\dots,0,-1,0,\dots,0)(12\dots m)(m+1,m+2,\dots,2m).\]
Then $y^m=(1,\dots,1,-1,\dots,-1)$, and $|y|=km$.
Furthermore,
\[y^b=(0,\dots,0,-1,0,\dots,0,1)(2m,\dots,m+1)(m\dots 1),\]
and
\[\begin{array}{rcl}
y^{-1}&=&(xa)^{-1}=a^{-1}x^{-1}\\
&=& (a^{-1}x^{-1}a)a^{-1}\\
&=&(0,\dots,0,-1,0,\dots,0,1)(2m,\dots,m+1)(m\dots 1)\\
&=& y^b.
%(m\dots 1)(2m,\dots,m+1)(-1,0,\dots,0,1,0,\dots,0),\\
%&=&(m\dots 1)(2m,\dots,m+1)(-1,0,\dots,0,1,0,\dots,0)(12\dots m)(m+1,m+2,\dots,2m)(2m,\dots,m+1)(m\dots 1)\\
\end{array}\]
Thus $y^b=y^{-1}$, and $\l y,b\r=\D_{2kn}$.
Since $y^m=(1,\dots,1,-1,\dots,-1)$ generates a semiregular subgroup on $\Sig$, the dihedral group $\l y,b\r$ is regular on $\Sig$, and so $X$ is a d-group.
\qed

\section{Proof of the main theorem}

By Theorem~\ref{qp-d-gps}, primitive d-groups are completely classified, and thus, to complete the proof of the main theorem, we only need to deal with the imprimitive case.
We first make a hypothesis.

\begin{hypothesis}\label{hypo-1}
{\rm
Let $G$ be a d-group of degree $2n$, with a regular dihedral subgroup $D$.
Let $\calB$ be a minimal block system and $B\in\calB$.
% such that $|B|>2$ and $|\calB|>2$.
Let $K=G_{(\calB)}$, the kernel of $G$ acting on $\calB$.
}
\end{hypothesis}

%Let $\calB$ be a $G$-invariant partition of $\Ome$, that is, $\calB$ is a block system for $G$.
%Then $\calB$ is also a block system for $D$ acting on $\Ome$, and $D^\calB$ is a transitive permutation group.
Since $D$ is a regular group on $\Ome$, the stabilizer $D_B$ induces a permutation group $D_B^B$, which is a regular group.
Regarding the action of $G$ on the block system $\calB$, the following lemma was given in \cite[Lemma 6.5]{et-Cay}.

\begin{lemma}\label{blocks-normal}
Let  $B\in\calB$, and let $K=G_{(\calB)}$ be the kernel of $G$ acting on $\calB$.
Then one of the following holds:
\begin{itemize}
\item[(i)] $D\cap K$ is regular on $B$, and $DK/K$ is regular on $\calB$, or

\item[(ii)] $D\cap K$ is semi-regular on $B$ with $2$ orbits and
$DK/K$ is transitive on $\calB$ with point stabilizer isomorphic to $\ZZ_2$.
In the latter case, $D/(D\cap K)$ has a cyclic regular subgroup.
\end{itemize}
In particular, $|D\cap K|=|B|$ or ${1\over2}|B|$.
\end{lemma}

Next, we further assume that $\calB$ is a minimal block system.
Then, for a block $B\in\calB$, the induced permutation group $G_B^B$ is primitive.

\begin{lemma}\label{D-on-B}
Using the notation defined above, the following hold:
\begin{itemize}
\item[(i)] $D_B^B\leq G_B^B$ is regular, and $G_B^B$ is a c-group or a d-group,

\item[(ii)] either $G_B^B$ is $2$-transitive, or $G_B^B=\ZZ_p\rtimes \ZZ_\ell<\AGL(1,p)$ with $|B|=p$ prime.

\end{itemize}
\end{lemma}
\proof
(i). Since $D$ is a regular group on $\Ome$ and $\calB$ is a block system for $D$, it follows that
$D_B^B\cong D_B$ is regular on $B$.
Since $D$ is dihedral, a subgroup $D_B$ is cyclic or dihedral.
Thus $G_B^B$ is a c-group or a d-group, respectively.

(ii). Since $\calB$ is a minimal block system, $G_B^B$ is a primitive permutation group by definition.
As $G_B^B$ is a c-group or a d-group by part~(i), by Theorems~\ref{qp-c-gps} and \ref{qp-d-gps}, we conclude that either $G_B^B$ is 2-transitive, or $|B|=p$ is a prime and $\ZZ_p\leq G_B^B<\AGL(1,p)$.
\qed

\begin{lemma}\label{K-2-gps}
If $K^B$ is an elementary abelian $p$-group with $p$, then so is $K$.
\end{lemma}
\proof
Since $K$ is normal in $G$, we have $K^B\lhd G_B^B$.
Label $\calB=\{B_1,B_2,\dots,B_m\}$.
Then $K\leq K^{B_1}\times K^{B_2}\times\dots\times K^{B_m}$, and $K^{B_1}\cong K^{B_2}\cong\dots\cong K^{B_m}$.
If $K^{B_i}$ is an elementary abelian $p$-group, then so is $K$.
\qed

\subsection{Small block case}

We handle the special case with $|B|=2$.
Let $Z$ be the cyclic subgroup of $D$ of index 2, and let $\s\in D\setminus Z$.
Then $|\s|=2$, and $D=Z\rtimes \l\s\r$.

\begin{lemma}\label{|B|=2}
Assume that $|B|=2$.
Then one of the following holds:
\begin{itemize}
\item[(i)] $K=1$, $|\calB|={1\over2}|\Omega|$, and $G^\calB\cong G$ is a c-group;
\item[(ii)] $K=\ZZ_2$, $G=\ZZ_2.G^\calB$, and either $D\cap K=1$ and $G^\calB$ is a c-group, and $D\cap K=\ZZ_2$ and $G^\calB$ is a d-group;
\item[(iii)] $K$ is an elementary abelian $2$-group of order at least $4$, and
there exists an orbital graph $\Ga_\calB[\ov\K_2]$, where $\Ga_\calB$ is an arc-transitive circulant or dihedrant.
\end{itemize}
Moreover, for each of the three cases, there indeed are examples.
\end{lemma}
\proof
Suppose $K\cap D=1$.
Then the action of $D$ on $\calB$ is faithful.
Thus the stabilizer of $B$ in $D$ is core-free in $D$, and so $D_B=\ZZ_2$.
Since $D_B$ is regular on $B$, we have $|B|=2$.
As $D^\calB$ is transitive, the cyclic subgroup $\l a\r$ is transitive and so regular on $\calB$.
Thus $G^\calB$ is a c-group.
In particular, if $K=1$, then $G^\calB$ is a c-group, as in part~(i).

Assume $K=\ZZ_2$.
Then $G=K.G^\calB=\ZZ_2.G^\calB$.
If $K\cap D=1$, then $G^\calB$ is a c-group as shown above;
if $K\cap D=\ZZ_2$, then $G^\calB$ is a d-group, as in part~(ii).

Now assume $|K|>2$.
Since $\Sym(B_i)=\S_2$ and $K\leq\Sym(B_1)\times\dots\times\Sym(B_m)=\S_2^m$, we conclude that $K$ is an elementary abelian 2-group.
Now the kernel of $K$ acting on a block $B$ is unfaithful, namely, $K_{(B)}\not=1$.
Thus $K_{(B)}$ is transitive on some $B'\in \calB$.
Pick $\a\in B$ and $\b\in B'$.
Let $\Ga$ be the orbital graph with arc set $(\a,\b)^G$.
Then $K_\a$ is transitive on $\Ga(\a)$, and it follows that $\Ga(\a)=B'$.
Thus we conclude that $\Ga=\Ga_\calB[\ov\K_2]$.

To complete the proof, we next construct examples of $d$-groups satisfying one of the conditions in the three cases.

{\bf Examples for part~(i):}
Let $G=\PGL(2,q)$, where $q=p^f$ with $p$ prime and $q\equiv3$ $(\mod 4)$.
Let $H=\ZZ_p^f\rtimes \ZZ_{q-1\over2}$, and $\Ome=[G:H]$.
Then $H$ is of odd order, and $|\Ome|=2(q+1)$.
Let $D=\D_{2(q+1)}$ be a maximal subgroup of $G$.
Then $H\cap D=1$ and $|G|=|H||D|$.
Hence $G=DH$, and $G$ is a d-group on $\Ome$ with $D$ being a dihedral regular subgroup.

Let $P$ be a maximal parabolic subgroup of $G$ which contains $H$, and let $\calB=[G:P]$.
Then $P=\ZZ_p^f\rtimes \ZZ_{q-1}$, and $\calB$ is block system for $G$ acting on $\Ome$.
As $|P|/|H|=2$, we conclude that $|\calB|={1\over2}|\Ome|$ and a block $B\in\calB$ has size 2.
Moreover, since the almost simple group $G$ is faithful on $\calB$, the d-group $G$ on $\Ome$ satisfies part~(i) with kernel $K=1$.

{\bf Examples for part~(ii):}
Let $G=\PGL(2,q)\times\l c\r$, where $q=p^f\equiv3$ $(\mod 4)$, and $|c|=\ZZ_2$.
Let $P=\ZZ_p^f\rtimes (\ZZ_{q-1\over2}\times\l z\r)$, a maximal parabolic subgroup of $G$, where $|z|=2$.
Let $H=\ZZ_p^f\rtimes (\ZZ_{q-1\over2}\times\l zc\r)\cong P$, and let $\Ome=[G:H]$.
Let $D=\D_{2(q+1)}$ be a maximal subgroup of $\PGL(2,q)$.
Then $H\cap D=1$ and $|G|=|H||D|$.
Hence $G=DH$, and $G$ is a d-group on $\Ome$ with $D$ being a dihedral regular subgroup.

Let $M=H\times\l c\r$, and let $\calB=[G:M]$.
Then, as $H$ is a subgroup of $M$ of index 2, $\calB$ is block system for $G$ acting on $\Ome$ and $|\calB|={1\over2}|\Ome|$.
Since the normal subgroup $\l c\r\lhd G$ is contained in the stabilizer $M$, we have $G_{(\calB)}=\l c\r$.
Thus the d-group $G$ on $\Ome$ satisfies part~(ii) with $K=\ZZ_2$ and $D\cap K=1$.

{\bf Examples for part~(iii)} are given in Construction~\ref{cons-1}.
This completes the proof of the lemma.
\qed

Finally, we present one more example for part~(ii) of Lemma~\ref{|B|=2}.

\begin{example}\label{K=2}
{\rm
Let $p\equiv3$ $(\mod 4)$ be a prime, and let $\Del=\{1,2,\dots,p\}$.
Let $G=\Sym(\Del)\times\l c\r$, where $|c|=2$.
Let $H=\Alt(\Del\setminus\{1\})$, and let $\Ome=[G:H]$.
Then $G$ is a permutation group on $\Ome$ of degree $4p$.
Let
\[D=\l(12\dots p),(2,p)(3,p-1)\dots({p+1\over2},{p+3\over2})\r\times\l c\r\cong\D_{2p}\times\ZZ_2=\D_{4p}.\]
Then $D$ is regular on $\Ome$, and $G$ is a d-group on $\Ome$.

Let $M=H\times\l c\r$, and let $\calB=[G:M]$.
Then, as $H$ is a subgroup of $M$ of index 2, $\calB$ is block system for $G$ acting on $\Ome$ and $|\calB|={1\over2}|\Ome|$.
Since the normal subgroup $\l c\r\lhd G$ is contained in the stabilizer $M$, we have $G_{(\calB)}=\l c\r$.
Thus the d-group $G$ on $\Ome$ satisfies part~(ii) with $K=D\cap K=\ZZ_2$.
}
\end{example}

\subsection{Bi-primitive d-groups}\

We here determine another extremal case -- bi-primitive d-groups.
A transitive permutation group $G\leq\Sym(\Ome)$ is said to be {\it bi-primitive} if there is a minimal block system with exactly two blocks, namely, $\calB=\{B_1,B_2\}$ such that $G_{B_i}$ is primitive on $B_i$.
Similar to the definition of quasiprimitive groups, a permutation group $G\leq\Sym(\Ome)$ is said to be {\it bi-quasiprimitive} if each minimal normal subgroup of $G$ has exactly two orbits.
Although the primitivity implies the quasiprimitivity, the bi-primitivity and the bi-quasiprimitivity are independent.
Bi-quasiprimitive d-groups have been classified in \cite{Song-Li}.
The following proposition determines bi-primitive d-groups.

Let $G$ be a bi-primitive permutation group on $\Omega$, and let $\calB$ be a block system.
Then $\calB$ contains exactly two blocks, namely, $\calB=\{B,B'\}$.
As usual, let $G^+=G_{B}=G_{B'}$.
Then $G^+$ is primitive on both $B$ and $B'$, and $G=G^+.\ZZ_2$.

\begin{proposition}\label{bi-primitive}
Assume that $G$ is a bi-primitive d-groups with notation defined above. Assume $R$ is a regular $d$-group containing in $G$.
Then one of the following holds.
\begin{itemize}
\item[(1)] there exists an orbital graph $\K_{m,m}$;

\item[(2)] $G=G^+\times\ZZ_2$, and either
\begin{itemize}
\item[(i)] $G^+=\S_n$ or $\A_{4m+1}$ with $n=4m+1$, and $R^+=\ZZ_n$, or
\item[(ii)] $G^+=\PGL(2,q).e$ and $R^+=\D_n$, where $n={q^d-1\over q-1}$ with $e\div f$ and $q=p^f$;
\end{itemize}

\item[(3)] $m=p$, $R=\D_{2p}$, $G=\ZZ_p\rtimes \ZZ_{2k}\le\AGL(1,p)$,
and $G_\o=\ZZ_k$, where $k$ is odd;

\item[(4)] $G=\S_4$ and $G^+=\A_4$;

\item[(5)] $G$ is almost simple, and $(G,R,G_\o)$ is one of the following triples.
\[\begin{array}{|llll|} \hline
G & R & G_\o & \mbox{condition} \\ \hline

\S_n &\D_{2n}& \A_{n-1} &\\

\PGL(2,q).e & \D_{2(q+1)} & \ZZ_r^f\rtimes \ZZ_{q-1\over 2}.\ZZ_e
 & e\div f\\

\PGL(2,q).2e & \D_{2(q+1)} & \ZZ_r^f\rtimes \ZZ_{q-1\over 2}.e &
 2e\div f, \mbox{$e$ odd} \\

\PGL(d,q).e.2 & \D_{2{q^d-1\over q-1}} & P_1 (\mbox{parabolic})
&d\geq 3,\ e\div f\\

\PGL(2,11) & \D_{22} & \A_5 & \\

\M_{12}.2 & \D_{24} & \M_{11} & \\ \hline
\end{array}\]
\end{itemize}
\end{proposition}
\proof
Assume first that $G_B$ is unfaithful on $B$.
Then $1\not=G_{(B)}\lhd G_B$, and hence $1\not=G_{(B)}^{B'}\lhd G_B^{B'}=G_{B'}^{B'}$.
Since $\calB$ is a minimal block system, $G_{B'}^{B'}$ is primitive, and so $G_{(B)}^{B'}$ is transitive.
It follows that there is an orbital graph $\K_{m,m}$, as in part~(1).

Suppose that $G^+=G_B\cong G_B^B$ is faithful on $B$.
Then $G=G^+.\ZZ_2$, and since $G_B^B$ is primitive by Theorem~\ref{mini-block}, either $|B|=p$ and $G^+\leq\AGL(1,p)$, or $G^+$ is a 2-transitive c-group or d-group.

{\bf Case 1.}\ First, assume $\C_G(G^+)=1$.
Then $G$ is isomorphic to a subgroup of $\Aut(G^+)$, and thus $G$ has a unique minimal normal subgroup, which is $\soc(G^+)$ and has exactly two orbits.
Thus $G$ is a bi-quasiprimitive $d$-group, and by \cite[Theorem~1.1]{Song-Li}, one of parts~(3)-(5) holds.

{\bf Case 2.}\ Next, assume $\C_G(G^+)\not=1$.
Since the center $\Z(G^+)$ is trivial, $G^+\cap \C_G(G^+)=1$.
As $G=G^+.\ZZ_2$, we have $\C_G(G^+)$ is of order 2, so $\C_G(G^+)=\l c\r=\ZZ_2$.
Then $G=G^+\times\l c\r$, where $|c|=2$, and $R^+:=R\cap G^+$ is a cyclic group or a dihedral group of index 2 in $R$.

Suppose that $R^+=\l a\r=\ZZ_n$ is cyclic.
Then $R=\l a\r \rtimes \l gc\r=\ZZ_n\rtimes \ZZ_2=\D_{2n}$ for some element $g\in G^+$.
Thus $a^{-1}=a^{gc}=a^g$, and so $\l a,g\r\cong\D_{2n}$ and $|g|=2$,
that is,   the c-group $G^+$ contains a transitive dihedral subgroup of order $2n$.
Analyzing the groups listed in Theorem~\ref{qp-c-gps}, we conclude that $G^+=\S_n$, $\A_{4m+1}$ with $n=4m+1$, or $\PGL(2,q).e$ with $n={q^d-1\over q-1}$, where $e\div f$ with $q=p^f$.
This is displayed in part~(2)(i).

Suppose that $R^+=\D_n$ is dihedral.
Then $R=\l R^+, gc\r=\D_{2n}$, where $g\in G^+$.
It follows that $gc$ has order $n$ or 2, and $\l R^+,g\r\cong\D_{2n}$.
Thus $|g|=1$, 2 or $n$.
If $|g|=1$, then $R=R^+\times \l c\r\cong\D_{2n}$, and so
$R^+=\D_n$ with $n/2$ odd.
Assume $|g|\not=1$.
Then there exists $x\in R^+$ such that $xg$ is an involution, and $\l R^+,xg\r=\D_{2n}$, namely, the d-group $G^+$ of degree $n$ has a transitive dihedral subgroup $\l R^+,xg\r$ of order $2n$.
It follows that the cyclic subgroup of $\l R^+,xg\r$ of order $n$ is regular, and so $G^+$ is also a c-group.
Arguing as in the previous paragraph, we have $G^+=\S_n$, $\A_{4m+1}$ with $n=4m+1$, or $\PGL(2,q).e$ with $n={q^d-1\over q-1}$, where $e\div f$ with $q=p^f$.
In this case, $R=\l R^+,xgc\r=R^+\rtimes\l xgc\r$???????
This is given in part~(2)(ii).
\qed

Examples of d-groups for part~(1) are classified in \cite{FanWW-1,FanWW-2}.
The d-groups satisfying part~(2) are given in the next example.

\begin{example}\label{bi-pri-d-1}
{\rm
\begin{itemize}
\item[(1)]
Let $\Delta$ be the set of 1-spaces of $\FF_q^2$, and let $G_0=\PGL(2,q)$.
Then $G_0$ is 2-transitive on $\Delta$ of degree $q+1$, and the point stabilizer is a parabolic subgroup $P=\AGL(1,q)$.
Let $D_0=\D_{2(q+1)}$ be a maximal subgroup of $H$.
Write $D_0=\ZZ_{q+1}\rtimes\l g\r$.
Let
\[G=G_0\times\l c\r,\ \Ome=[G:P],\ D=D_0\rtimes\l gc\r.\]
Then $G$ is bi-primitive on $\Ome$, and $D$ is a regular dihedral subgroup of $G$ on $\Ome$,
namely, $G=\PGL(2,q)\times\ZZ_2$ is a bi-primitive d-group of degree $2(q+1)$.

\item[(b)] Let $G_0=\A_n=\A_{4m+1}$, acting on $\Delta=\{1,2,\dots,4m+1\}$ naturally.
Let $a=(12\dots 4m+1)$, and $g=(1,4m)(2,4m-1)\dots(2m,2m+1)$.
Then $\l a\r$ is regular on $\Delta$, and $\l a,g\r=\D_{2n}$ is transitive on $\Delta$.
Let
\[G=G_0\times\l c\r,\ \Ome=[G:\A_{4m}],\ D=D_0\rtimes\l gc\r.\]
Then $G$ is bi-primitive on $\Ome$, and $D$ is a regular dihedral subgroup of $G$ on $\Ome$,
namely, $G=\A_{4m+1}\times\ZZ_2$ is a bi-primitive d-group of degree $2(4m+1)$.

\item[(c)] Similarly, the symmetric group $\S_n$ acting on $n$ points gives rise to a bi-primitive d-group of $2n$.
\end{itemize}
}
\end{example}

Examples of d-groups satisfying part~(3) are all as follows.

\begin{example}\label{c-d-gps}
{\rm
Let $G=\l a\r\rtimes \l b\r=\ZZ_p\rtimes \ZZ_\ell\le\AGL(1,p)$, where $p$ is an odd prime and $\ell$ is an even divisor of $p-1$.
Let $H=\l b^2\r=\ZZ_{\ell/2}$, and let $\Ome=[G:H]$.
Then $R=\l a\r\rtimes \l b^{\ell/2}\r=\D_{2p}$ is regular on $\Ome$, and $G$ is a d-group on $\Ome$.
Let $M=\l b\r$.
Then $H<M$, and $\calB=[G:M]$ is a $G$-invaraint partiton of $\Ome$.
Now $Z=\l a\r=\ZZ_p$ is regular on $\calB$, and $G^\calB\cong G$ is a c-group on $\calB$.
}
\end{example}

\subsection{General case}\

In this subsection, we handle the general case.
Let $\calB$ be a minimal block system for $G$ acting on $\Ome$, satisfying Hypothesis~\ref{hypo-1} with $|B|>2$ and $|\calB|>2$.

%\begin{hypothesis}\label{hypo}
%{\rm
%Let $G$ be a d-group of degree $2n$, with a regular dihedral subgroup $D$.
%Let $\calB$ be a minimal block system and $B\in\calB$ such that $|B|>2$ and $|\calB|>2$.
%Let $K=G_{(\calB)}$, the kernel of $G$ acting on $\calB$.
%}
%\end{hypothesis}

\begin{lemma}\label{K-not=1}
The permutation group $K^B$ is a transitive subgroup of $G_B^B$, and $D\cap K$ is cyclic.
 %of order $|B|$ or ${1\over2}|B|$.
\end{lemma}
\proof
By Lemma~\ref{D-on-B}, as $|B|>2$, we have $|R\cap K|>1$, and so $K\not=1$.
Thus $1\not=K^B\lhd G_B^B$.
Since $G_B^B$ is primitive, the normal subgroup $K^B$ is transitive.

Since $R\cap K\lhd R$ and $R/(R\cap K)$ had order divisible by $|\calB|\geq 3$, we conclude that $R\cap K$ is cyclic.
\qed

The following lemma characterizes the kernel $K=G_{(\calB)}$.

\begin{lemma}\label{K-primitive}
One of the following holds:
\begin{itemize}
\item[(i)] $K^B$ is $2$-transitive,

\item[(ii)] $|B|=p$ with $p$ prime, and $\ZZ_p\lhd K^B\lhd G_B^B=\ZZ_p\rtimes \ZZ_\ell<\AGL(1,p)$,

\item[(iii)] $|B|=4$ and $K$ is an elementary abelian $2$-group.
\end{itemize}
\end{lemma}
\proof
By Lemma~\ref{D-on-B}, the induced permutation group $G_B^B$ is a c-group or a d-group, and
by Lemma~\ref{K-not=1}, $K^B$ is a transitive normal subgroup of $G_B^B$.
If $K^B$ is primitive, then by Theorems~\ref{qp-c-gps} and \ref{qp-d-gps}, either $K^B$ is 2-transitive, or $\ZZ_p\lhd K^B\lhd G_B^B=\ZZ_p\rtimes \ZZ_\ell<\AGL(1,p)$, as in part~(i) or (ii).

Now assume that $K^B$ is imprimitive.
Then $G_B^B$ is affine and the socle $\soc(G_B^B)$ is not simple.
By Theorem~\ref{qp-d-gps}, we have $\soc(G_B^B)=\ZZ_2^2$, or $\ZZ_2^3$, or $\ZZ_2^4$.

Suppose that $\soc(G_B^B)=\ZZ_2^3$ or $\ZZ_2^4$.
Then the block size $|B|=8$ or 16, respectively.
On the other hand, as $|\calB|>2$, the normal subgroup $R\cap K$ of the dihedral group $R$ is of index bigger than 2, and so $R\cap K$ is cyclic, say $K\cap R=\ZZ_\ell$, with $\ell\ge |B|/2\ge4$.
Since $R$ is regular on $\Ome$, it follows that $K\cap R$  is faithful on $B$, and
$K^B\lhd G_B^B$ contains cyclic subgroup $(K\cap R)^B$ of order 4.
By Theorem~\ref{qp-d-gps}, the primitive permutation group $G_B^B=N.L$, where $N=\soc(G_B^B)$ and $L$ is almost simple.
Thus, $K^B=N.L_0$ is primitive where $L_0\ge\soc(L)$, which is a contradiction.

We therefore have $|B|=4$, and $G_B^B=\A_4$ or $\S_4$.
It follows that $K$ is a $\{2,3\}$-group.
If $|K|$ is divisible by 3, then $K^B$ is of order divisible by 3, and
as $K^B\lhd G_B^B$, we conclude that $K^B=\A_4$ or $\S_4$, and thus $K^B$ is primitive.
This contradiction shows that $K$ is a 2-group.
Suppose that $K$ has an element $g$ of order 4.
Then $g^2$ acts non-trivially on some block, say $B$, and thus $g$ is an element of $K^B$
of order 4.
Noticing that $K^B$ is a normal subgroup of the primitive permutation group $G_B^B$, we conclude that $K^B=G_B^B=\S_4$, which is a contradiction.
Hence, every element of $K$ is of order 2, and $K$ is an elementary abelian 2-group.
\qed

We analyse the three cases separately.

\begin{lemma}\label{K-2-trans}
Assume that $K^B$ is primitive.
Then one of the following holds:
\begin{itemize}
\item[(1)] $K$ is unfaithful on $B$, and there  exists an orbital graph for $G$ of the form $\Ga_\calB[\ov\K_b]$, where $b=|B|$;

\item[(2)] $\ZZ_p\rtimes \ZZ_k=K\cong K^B\leq\AGL(1,p)$, where $k$ is an odd divisor of $p-1$, and either
\begin{itemize}
\item[(i)] $G=((\ZZ_p\rtimes \ZZ_k)\times H).\ZZ_{2\ell}$, or

\item[(ii)] $K=\ZZ_p$, and $G=H.\ZZ_{2\ell}$ where $H=\C_G(K)$ is a non-split central extension of $\ZZ_p$ by $H^\calB$,

\end{itemize}
where $\ell$ is an odd divisor of $(p-1)/k$ and $\ZZ_k.\ZZ_{2\ell}=\ZZ_{2k\ell}\leq\ZZ_{p-1}$.

\item[(3)] $K\cong K^B$ is non-solvable, and $G=(K\times H).\calO$, where $\ZZ_2\leq\calO\leq\Out(K)$;

\end{itemize}
\end{lemma}
\proof
If $K$ is unfaithful on $B$, then, by Lemma~\ref{key-lem}, there  exists an orbital graph for $G$ of the form $\Ga_\calB[\ov\K_b]$, where $b=|B|$, as in part~(1).
Assume that $K$ is faithful on $B$ in the following.
Since $|\calB|,|B|>2$, it follows from Lemma~\ref{blocks-normal} that the intersection $K\cap D\not=1$.

{\bf Case 1.}\ \ First, assume that $\ZZ_p\leq K\cong K^B\leq\AGL(1,p)$.
Then $K=\ZZ_p\rtimes \ZZ_k$ with some integer $k\div(p-1)$.
%As $D\cap K\not=1$, it follows that $K=\ZZ_p\rtimes \ZZ_k$ with $k$ odd.
Let $H=\C_G(K)$, the centralizer of $K$ in $G$.
Since $K\lhd G$, we have $H,KH\lhd G$, and hence $G/(KH)\leq\Out(K)$.

Suppose that $k\not=1$.
Then $K=\ZZ_p\rtimes \ZZ_k$ is a Frobenius group, and so $K\cap H=1$.
Thus $\Aut(K)=\ZZ_p\rtimes\ZZ_{p-1}=\AGL(1,p)$, and $\Out(K)=\ZZ_{p-1}/\ZZ_k=\ZZ_{(p-1)/k}$.
So $H\times K=HK\lhd G$, and
\[G=(H\times K).\ZZ_\ell,\]
where $\ell\div (p-1)/k$.
Since $|\calB|>2$, there are at least 3 orbits of $K\cap D$ on $\Ome$.
Thus $K\cap D$ is cyclic and normal in $D$.
It follows that $K\cap D=\ZZ_p$, and $H=\C_G(K\cap D)$.
Thus $R\cap (K\times H)$ is the cyclic subgroup of $D$ of index 2, and
\[G=(K\times H).\ZZ_{2\ell}=((\ZZ_p\rtimes\ZZ_k)\times K).\ZZ_{2\ell}.\]
for some integer $\ell$, as in part~(2)(i).

Now suppose that $k=1$ and $K=\ZZ_p$.
If $K.H^\calB$ is a split extension, then $KH=K\times H^\calB$, as in part~(2)(i) with $k=1$.
If the extension $K.H^\calB$ is non-split, then the group $G$ satisfies part~(2)(ii).

{\bf Case 2.}\ \ Next, assume that $K\cong K^B$ is non-solvable.
Let $H=\C_G(K)$, the centralizer of $K$ in $G$.
Since $K\cong K^B$ is 2-transitive, we have $K\cap H=1$, and hence $KH=K\times H\lhd G$.
The factor group $G/(KH)\leq\Out(K)$.
Suppose that $D\leq K\times H$.
Then $D\leq D_1\times D_2$, where $D_1$ is the projection of $D$ in $K$ and $D_2$ is the projection of $D$ in $H$.
Since $D_1$ and $D_2$ are dihedral groups of order at least 4, this is not possible.
Thus $D$ is not contained in $K\times H$, and $\calO\geq\ZZ_2$.
In particular, $\Out(K)\geq\ZZ_2$.
By Theorems~\ref{qp-c-gps} and \ref{qp-d-gps}, either $K=\A_4$ or $2^4\rtimes \A_6$, or $K$ is one of the following almost simple groups:
\[\A_{2m+1},\ \M_{22},\ \A_{4m},\ \PSL(2,q).e\ \mbox{with $q=r^f\equiv3$ $(\mod 4)$ and $e\div f$.}\]
Moreover, since $K\cap D$ is cyclic, we conclude that $K=\A_{2m+1}$, $\M_{22}$ or $\A_{4m}$.
\qed

\begin{lemma}\label{mini-actions}
If $|B|=4$ and $K$ is an elementary abelian group, then $R\cap K=\ZZ_2$, $G=2^d.H$ where $H$ is a c-group.
\end{lemma}
\proof
Since $|\calB|\geq4$, the intersection $R\cap K$ is cyclic, and so $R\cap K\cong\ZZ_2$.
By Lemma~\ref{blocks-normal}, the induced permutation group $G^\calB\cong G/K$ is a c-group of degree $n/4$.
\qed

We end this subsection with some examples for the case where $|B|=4$ and $K$ is an elementary abelian 2-group,
$K\cap R=\ZZ_2$, and $G^\calB$ is a c-group and contains a transitive dihedral subgroup.

\begin{example}\label{Alt6}
{\rm
Let $G=\ZZ_2^9\rtimes \PGammaL(2,9)$, and let $H=\ZZ_2^7\rtimes 3^2\rtimes [2^4]$.
Let $\Ome=[G:H]$, degree 40.
Let $M=\ZZ_2^9\rtimes 3^2\rtimes [2^4]>H$.
Then $\calB:=[G:M]$ forms a $G$-invariant partiton of $\Ome$, and
$G^\calB=\PGammaL(2,9)$, of degree 10, and $G_B^\calB=3^2\rtimes [2^4]$.
Now $G^\calB$ has a cyclic regular subgroup $Z^\calB\cong\ZZ_{10}$, and
a transitive dihedral subgroup $R^\calB=\D_{20}$.
}
\end{example}

\begin{example}\label{7-6}
{\rm
Let $G=2^6\rtimes (7\rtimes 6)$ is an irreducible subgroup of $\AGL(6,2)$.
Let $G_\o=2^4\rtimes 6$, and let $K=\soc(G)=2^6$.
Then $\Ome:=[G:G_\o]$ is of degree $28$, and $G$ has a dihedral subgroup $R\cong\D_{28}$.
Let $\calB$ be the set of $K$-orbits on $\Ome$.
Then $|\calB|=7$, and $B\in\calB$ has size 4.
Furthermore, $G_B=2^6\rtimes 6$, and $G_B^B=\A_4$, and $G^\calB=7\rtimes 6$.
}
\end{example}

\subsection{Proof of Theorem~\ref{thm-1}}\

We here summarize what we have obtained to complete the proof of Theorem~\ref{thm-1}.

Let $G$ be an imprimitive d-group on $\Ome$ of degree $2n$.
Let $\calB$ be a minimal block system for $G$ on $\Ome$, and let $K=G_{(\calB)}$.
We will list the candidates for $G$ according to the property of $K$.

If $K=1$, then $|B|=2$ by Lemma~\ref{K-not=1}, and $G\cong G^\calB$ is also a c-group by Lemma~\ref{|B|=2}, as in part~(1).

Assume $K\not=1$.
If $K^B$ is imprimitive, then by Lemmas~\ref{K-primitive} and \ref{mini-actions}, we have that $|B|=4$, and $K$ is an elementary abelian 2-group, and $G^\calB$ is a c-group of degree $n/4$, as in part~(5).
Suppose that $K^B$ is primitive.
If $K_{(B)}\not=1$, then there exists an orbital graph $\Ga=\Ga_\calB[\ov\K_b]$ where $b=|B|$, as in part~(2).

Finally, suppose that $K\cong K^B$ is faithful and primitive.
If $K=\ZZ_p$ with $p$ prime, then by Lemmas~\ref{|B|=2} and \ref{K-2-trans}, part~(3) of Theorem~\ref{thm-1} occurs.
On the other hand, the centre $\Z(K)=1$ for $K\not=\ZZ_p$ by Lemma~\ref{K-primitive}.
In this case, by Lemma~\ref{K-2-trans}, part~(4) of Theorem~\ref{thm-1} occurs.
\qed

\end{document}